\documentclass[12pt]{amsart}

\usepackage{verbatim, amscd,epsfig,amsmath,amsthm,amstext,amssymb,a4}
\usepackage{enumerate}
\usepackage{mathtools}
\usepackage[all]{xy}
\usepackage[hyperindex=true,plainpages=false,colorlinks=false,pdfpagelabels]{hyperref}
\usepackage{verbatim}
\usepackage{amssymb}

\numberwithin{equation}{section}

\newtheorem{theorem}{Theorem}[section]
\newtheorem{proposition}[theorem]{Proposition}
\newtheorem{lemma}[theorem]{Lemma}

\newtheorem{question}[theorem]{Question}

\theoremstyle{definition}

\newtheorem{remark}[theorem]{Remark}

\newcommand{\R}{\mathbb{R}}
\newcommand{\C}{\mathbb{C}}
\usepackage{graphicx,xcolor}

\begin{document}

\title[Holomorphic Higgs bundles on Teichm\"uller space]{Holomorphic
Higgs bundles over the Teichm\"uller space}

\author[I. Biswas]{Indranil Biswas}

\address{Department of Mathematics, Shiv Nadar University, NH91, Tehsil Dadri,
Greater Noida, Uttar Pradesh 201314, India}

\email{indranil.biswas@snu.edu.in, indranil29@gmail.com}

\author[L. Heller]{Lynn Heller}

\address{HIMIS, Shenzhen, China}

\email{lynnheller@himis-sz.cn}

\author[S. Heller]{Sebastian Heller}

\address{Department of Mathematics and Institute of Mathematical Sciences,
The Chinese University of Hong Kong}

\email{s.g.heller@icloud.com}



\begin{abstract}
We study those representations $\rho$ of the fundamental group of a compact connected oriented
surface $X$ such that the Higgs data corresponding to $\rho$ --- via the non-abelian Hodge
correspondence --- have the property that they depend holomorphically on the Riemann surface
structure $\Sigma\,=\, (X,\, J)$ on $X$. For representations $\rho$ into $\mathrm{SL}(2,
\mathbb C)$ we show that this holomorphic dependency is equivalent to $\rho$ being
unitary. For higher ranks this equivalence fails --- we 
show the existence of non-unitary and irreducible representations of the
fundamental group of a compact connected oriented
surface into $\mathrm{SL}(n,\mathbb C)$ admitting 
Higgs data that are holomorphic in $\Sigma$, for $n$ large enough.
\end{abstract}

\maketitle

\section{Introduction}

The celebrated non-abelian Hodge correspondence gives a real analytic diffeomorphism between the moduli space 
of stable Higgs pairs of rank $n$
on a given compact connected Riemann surface $\Sigma \,= (X,\, J)$ and the moduli space of 
irreducible flat connections of rank $n$ on $\Sigma$. The latter can be identified with the Betti moduli
space of the underlying smooth surface $X$ for the group
$\text{GL}(n,{\mathbb C})$, i.e., the moduli space of representations of $\pi_1(X,\,x_0)$
into $\text{GL}(n,{\mathbb C})$. Note that the 
Higgs bundle moduli space depends on the Riemann surface structure (in fact there is a Torelli-type theorem 
asserting that the isomorphism class of the moduli space actually determines the
isomorphism class of the Riemann surface) while the Betti moduli space is independent 
of the complex structure on the surface.\\

A natural question is therefore how the non-abelian Hodge correspondence depends on the Riemann surface 
structure. In this paper, we address the following special case: Given a representation $\rho$ of the fundamental group 
of $X$, when does the corresponding Higgs pair depend holomorphically on the Riemann surface structure 
$\Sigma$. We first observe in Proposition \ref{pro:un} that this property clearly holds for
all unitary representations of $\pi_1(X,\,x_0)$.
The remainder of Section \ref{sec:set} presents the relevant notation.

In Section \ref{sec:rk2} we prove that 
the only representations $\rho$ of $\pi_1(X,\,x_0)$
to $\mathrm{SL}(2,\mathbb C)$ for which the Higgs data are holomorphic in 
$\Sigma$ are the unitary representations (see Theorem \ref{The2}). The proof of this theorem makes use of the 
energy functional of the associated harmonic map, which is defined over the Teichm\"uller space.\\

In the last section, we provide examples of non-unitary and irreducible representations of
$\pi_1(X,\,x_0)$ into 
$\mathrm{SL}(n,\mathbb C)$ such that the Higgs data depend holomorphically on the complex structure of the 
surface, when $n$ is large enough (see Theorem \ref{mthm}). Some natural questions
that arise in these investigations are as follows (see also \cite{DT}):

\begin{question}
How to classify representations of $\pi_1(X,\,x_0)$ with the property that the Higgs data 
depend holomorphically on the complex structure? Do these representations admit any (branched) minimal 
surfaces, or must the associated energy over Teichm\"uller space even be constant?
\end{question}

\section{The set-up}\label{sec:set}

Let $X$ be a compact connected oriented $C^\infty$-surface of genus $g$, with $g\, \geq\, 2$. Let ${\rm 
Diff}^+(X)$ denote the group of all orientation preserving diffeomorphisms of $X$ and ${\rm 
Diff}^0(X)\,\subset\,{\rm Diff}^+(X)$ the normal subgroup consisting of all diffeomorphisms homotopic to the 
identity map of $X$. The quotient
\begin{equation}\label{e1}
{\rm MCG}(X)\,\,:=\,\, {\rm Diff}^+(X)/{\rm Diff}^0(X)
\end{equation}
is known as the mapping class group of $X$.\\

Let $n\,\in\, \mathbb N_{\geq2}$ be an integer, and consider the complex Lie group $G$, which we will fix in this 
paper to be either $\text{SL}(n, {\mathbb C})$ or $\text{PGL}(n, {\mathbb C})$. Fix a base point $x_0\, 
\in\, X$, and consider a homomorphism $$\rho\ :\ \pi_1(X,\, x_0)\ \longrightarrow\ G$$
of the fundamental group. By composing the 
standard action of $G$ on ${\mathbb C}{\mathbb P}^{n-1}$ with $\rho$ we obtain an action of the fundamental group 
$\pi_1(X,\, x_0)$ on ${\mathbb C}{\mathbb P}^{n-1}$. The homomorphism $\rho$ is called irreducible if no proper 
linear subspace of ${\mathbb C}{\mathbb P}^{n-1}$ is preserved by this action of $\pi_1(X,\, x_0)$ on
${\mathbb C}{\mathbb P}^{n-1}$. The 
homomorphism $\rho$ is called completely reducible if the Zariski closure of ${\rm image}(\rho)$ is a 
reductive subgroup of $G$. In particular, every irreducible homomorphism is completely reducible. \\

denote by $\text{Hom}(\pi_1(X,\, x_0), \, G)$ the space of homomorphisms from $\pi_1(X,\, x_0)$ to $G$.
Let $$\text{Hom}'(\pi_1(X,\, x_0), \, G)\ \subset\ \text{Hom}(\pi_1(X,\, x_0), \, G)$$ be the space of
completely reducible homomorphisms from $\pi_1(X,\, x_0)$ to $G$. The conjugation action of $G$ on itself
gives rise to an action of $G$ on $\text{Hom}'(\pi_1(X,\, x_0), \, G)$ and the corresponding quotient space
\begin{equation}\label{e2}
{\mathcal R}_X(G) \ :=\ \text{Hom}'(\pi_1(X,\, x_0), \, G)/G
\end{equation}
is known as the $G$--character variety of $\pi_1(X,\, x_0)$. Via monodromy representation, the space
${\mathcal R}_X(G)$
can be identified with the isomorphism classes of flat $G$--connections $(E,\, \nabla)$ over
$X$ satisfying the condition that the monodromy of $(E,\, \nabla)$ generates a reductive subgroup, meaning the
closure of the image of the monodromy homomorphism is a reductive subgroup of $G$.

Let $\text{Com}(X)$ denote the space of all $C^\infty$-complex structures on $X$ compatible
with its orientation; it is equipped with a natural complex structure that is preserved by the natural action of ${\rm Diff}^+(X)$
on $\text{Com}(X)$. The quotient
\begin{equation}\label{dx}
{\mathcal T}(X)\ :=\ \text{Com}(X)/{\rm Diff}^0(X)
\end{equation}
is the Teichm\"uller space associated to $X$, which is a contractible complex manifold of
dimension $6g-6$. The action of ${\rm Diff}^+(X)$ on $\text{Com}(X)$ induces
an action of the mapping class group ${\rm MCG}(X)$ (see \eqref{e1}) on ${\mathcal T}(X)$
via holomorphic automorphisms.\\

Now let
${\rm Diff}^+_1(X)$ be the group of orientation preserving diffeomorphisms $f$ of $X$ fixing the base point,
i.e., $f(x_0)\,=\, x_0$, and 
let
$$
{\rm Diff}^0_1(X)\,\,\subset\,\, {\rm Diff}^+_1(X)
$$
be the connected component of ${\rm Diff}^+_1(X)$ containing the identity element, which is clearly a
normal subgroup of ${\rm Diff}^+_1(X)$. Furthermore, the natural action of ${\rm Diff}^+_1(X)$ on
$\text{Com}(X)$ preserves its complex structure, and the
quotient $$\widehat{Z}\ :=\ \text{Com}(X)/{\rm Diff}^0_1(X),$$ which is a contractible
complex manifold of complex dimension $6g-5$, admits a natural holomorphic projection map to the
Teichm\"uller space 
\begin{equation}\label{e6a}
\widehat{\phi} \,\,\colon\,\, \widehat{Z}\,:=\, \text{Com}(X)/{\rm Diff}^0_1(X)\,
\longrightarrow\,\, {\mathcal T}(X)\,=\, \text{Com}(X)/{\rm Diff}^0(X).
\end{equation}
See \cite{FM,Bi} for more details. Define
\begin{equation}\label{e7}
{\rm MCG}^1(X)\,\,:=\,\, {\rm Diff}^+_1(X)/{\rm Diff}^0_1(X).
\end{equation}
From the definition of ${\rm MCG}^1(X)$ in \eqref{e7} it follows immediately that the action
of ${\rm Diff}^+_1(X)$ on $\text{Com}(X)$ produces an action of ${\rm MCG}^1(X)$ on
$\widehat Z$ (see \eqref{e6a}) via holomorphic automorphisms.
Moreover, we have the following short exact sequence of groups
\begin{equation}\label{e8}
e\, \longrightarrow\, \pi_1(X,\, x_0) \, \longrightarrow\, {\rm MCG}^1(X)
\, \stackrel{\Pi}{\longrightarrow}\, {\rm MCG}(X) \, \longrightarrow\, e,
\end{equation}
which is known as the Birman exact sequence (see \cite[pp.~213--238]{Bi},
\cite[pp.~97--98, Theorem 4.6]{FM}). For any $z\, \in\, \widehat{Z}$ and $\gamma\, \in\, {\rm MCG}^1(X)$,
we have
\begin{equation}\label{f3}
\widehat{\phi}(\gamma\cdot z) \,\,=\,\, \Pi (\gamma)\cdot \widehat{\phi}(z),
\end{equation}
where $\widehat\phi$ and $\Pi$ are the maps constructed in \eqref{e6} and \eqref{e8} respectively.
In other words, the map $\widehat\phi$ is equivariant with respect to the homomorphism $\Pi$ of groups.\\

As noted above, ${\rm MCG}^1(X)$ acts on $\widehat Z$
in \eqref{e6a} via holomorphic automorphisms. Define
$$
Z \ \, :=\ \, \widehat{Z}/\pi_1(X,\, x_0),
$$
where $\pi_1(X,\, x_0) \, \subset\, {\rm MCG}^1(X)$ is the subgroup in \eqref{e8}.
The projection $\widehat{\phi}$ in \eqref{e6a} clearly factors through the quotient map
\begin{equation}\label{uc0}
\eta\ :\ \widehat{Z} \ \longrightarrow\ \widehat{Z}/\pi_1(X,\, x_0)\ =\ Z.
\end{equation}
Consequently, $\widehat\phi$ induces a natural holomorphic projection map to the Teichm\"uller space 
\begin{equation}\label{e6}
\phi \,\,\colon\,\, Z\,\,:=\,\, \widehat{Z}/\pi_1(X,\, x_0)\,\,
\longrightarrow\,\, {\mathcal T}(X);
\end{equation}
in other words, $\phi\circ\eta\,=\,\widehat\phi$. Note that the map $\eta$ in
\eqref{uc0} is a universal covering (recall that $\widehat Z$ is contractible).

The fibers of $\phi$ in \eqref{e6} are compact Riemann surfaces of genus $g\,=\, \text{genus}(X)$,
and $Z$ is a complex manifold of dimension
$6g-5$; it is in fact the universal Riemann surface over ${\mathcal T}(X)$.
The fiber of $\phi$ over any point $\Sigma\, \in\, {\mathcal T}(X)$
has a natural biholomorphism with the compact Riemann surface $\Sigma$.

For any $z_0\, \in\, Z$, the homomorphism $$\pi_1(\phi^{-1}(\phi(z_0)),\, z_0)\ \longrightarrow\
\pi_1(Z,\, z_0)$$ given by the inclusion map $\phi^{-1}(\phi(z_0))\, \hookrightarrow\, Z$ is an
isomorphism; this is because ${\mathcal T}(X)$ is contractible. Furthermore, $\pi_1(\phi^{-1}
(\phi(z_0)),\, z_0)$ can be identified with $\pi_1(X,\, x_0)$ up to inner automorphisms.
Consequently, any completely reducible homomorphism
\begin{equation}\label{er}
\rho\, :\, \pi_1(X,\, x_0) \, \longrightarrow\, G
\end{equation}
produces a flat $G$--connection on $Z,$ which we denote by 
\begin{equation}\label{e9}
(E_\rho,\, \nabla^\rho).
\end{equation}
For any $z\, \in\, Z$, consider the restriction of the flat bundle in \eqref{e9} to the fiber
$\phi^{-1}(z)$. The non-abelian Hodge correspondence associates to it a $G$--Higgs bundle on
$\phi^{-1}(z)$. In this way, we obtain a family of Higgs bundles parametrized by ${\mathcal T}(X)$.

Let $${\mathcal H}_G\,\, \longrightarrow\,\, {\mathcal T}(X)$$
denote the moduli space of semistable $G$--Higgs bundles over ${\mathcal T}(X)$. In other
words, ${\mathcal H}_G$ parametrizes all pairs $(Y,\, (E,\,\Phi))$, where $Y\, \in\, \mathcal T(X)$ and
$(E,\, \Phi)$ is a semistable $G$--Higgs bundle on $Y$. This ${\mathcal H}_G$ is a complex analytic space.
This structure of complex analytic space is uniquely determined by the following condition:

Let $\varpi\, :\, Y_S\, \longrightarrow\, S$
be a holomorphic family of compact Riemann surfaces of genus $g$ parametrized by a complex manifold $S$ such
that for every point $s_0\, \in\, S$, the fiber $\varpi^{-1}(s_0)$ is identified with the surface $X$ in 
\eqref{dx} by a class of orientation preserving diffeomorphism that differ by ${\rm Diff}^0(X)$; this condition means that
there is a holomorphic map $\widetilde{\varpi}\, :\, S\, \longrightarrow\, {\mathcal T}(X)$, and $Y_S$ is the pulled-back family of Riemann surfaces $\widetilde{\varpi}^*Z$, where $Z$ is the holomorphic family of Riemann surfaces in
\eqref{e6}. Let $({\mathcal E}_S,\, \vartheta_S)$ be a holomorphic family of semistable $G$--Higgs bundles on $X_S$.
Then the classifying map $S\, \longrightarrow\, {\mathcal H}_G$ is holomorphic.

The underlying real space for ${\mathcal H}_G$ is ${\mathcal T}\times {\rm Hom}(\pi_1(X,\, x_0),\, G)/\!\!/G$. The
complex structure on it is uniquely determined by the following two conditions:
\begin{enumerate}
\item The projection ${\mathcal H}_G\,\, \longrightarrow\,\, {\mathcal T}(X)$ is holomorphic.

\item If $Y_S\, \longrightarrow\, S$ is a holomorphic family of marked compact Riemann surfaces of genus
$g$, and $({\mathcal E},\, \vartheta)$ is a relative holomorphic $G$--Higgs bundle on $Y_S$ which is
fiberwise semistable, then the classifying map $S\, \longrightarrow\, {\mathcal H}_G$ is holomorphic.
\end{enumerate}
Let
\begin{equation}\label{eP}
P\, :\, {\mathcal H}_G\, \longrightarrow\, \mathcal T(X),\, \ \ (Y,\, (E,\,\Phi))\, \longmapsto\, Y
\end{equation}
be the natural holomorphic projection. For any $\rho$ as in \eqref{er}, we consider the section of
$(\mathcal H_G, P)$
\begin{equation}\label{eP2}
{\mathbb H}_\rho \ \colon \ \mathcal T(X) \ \longrightarrow \ {\mathcal H}_G
\end{equation}
that assigns to every surface $\Sigma\,\in\, \mathcal T(X)$ the Higgs bundle over
$\Sigma$ corresponding to the flat $G$--connection $(E_\rho, \,\nabla_\rho)\big\vert_{\phi^{-1} (\Sigma)}$
given in \eqref{e9}. This section is $C^\infty$; see \cite[p.~7823, \S~1.2.2]{Tos}
(see also \cite[p.~7834, Theorem 4.7]{Tos} and \cite[Section 4.4.2]{Tos}).

The question we want to address is the following: For which representations $\rho$
the section ${\mathbb H}_\rho$ in \eqref{eP2} turns out to be holomorphic?

See \cite[p.~893, Theorem 3.1]{B97} for a result similar to Proposition \ref{pro:un} for the special
case of unitary representations.

\begin{proposition}\label{pro:un}
For any unitary $\rho$ the section ${\mathbb H}_\rho$ in \eqref{eP2} is holomorphic.
\end{proposition}

\begin{proof}
Let $\rho$ be a unitary representation. For any $\Sigma\, \in\, T(X)$, consider ${\mathbb H}_\rho(\Sigma)
\,=\, (\Sigma,\, (E_\rho,\,\Phi_\rho))$. Then $\Phi_\rho\,=\, 0$ by unitarity, and $E_\rho$ is the holomorphic
principal $G$--bundle on $\Sigma$ given by $\rho$. Consequently, the section
${\mathbb H}_\rho$ is holomorphic. 
\end{proof}

For the case of the representations $\rho$ to $\mathrm{SL}(2,\mathbb C)$, we show in Lemma \ref{The1} and 
Theorem \ref{The2} that ${\mathbb H}_\rho$ is holomorphic if and only if $\rho$ is unitary. The proof uses 
properties of the energy functional $\mathcal E_\rho$ on the Teichm\"uller space determined by the 
Dirichlet energy of the associated harmonic maps into the corresponding symmetric space; this symmetric 
space coincides with the hyperbolic 3-space
\begin{equation}\label{h3s}
H^3\ =\ \mathrm{SL}(2,{\mathbb C})/\mathrm{SU}(2).
\end{equation}
In particular, we have the following: The energy functional $\mathcal E_\rho$, which is bounded from below by
$0$, cannot be proper if ${\mathbb H}_\rho$ is holomorphic.

However, Theorem \ref{The2} does not generalize
to arbitrary rank. Indeed, it is shown in Section \ref{sec:hr} that there are 
non-unitary completely reducible representations $\rho\,\colon\,\pi_1(X,\,x_0)\,\
\longrightarrow\, \mathrm{SL}(n,\mathbb C)$ --- for some large 
$n$ --- satisfying the condition that the section ${\mathbb H}_\rho$ is holomorphic.

\section{Higgs Bundles of rank two}\label{sec:rk2}

In this section we restrict ourselves to $G \,= \,{\rm SL}(2, \C)$. We first start with an introductory 
example, where the family of Higgs bundles is never holomorphic, to illustrate our arguments.

We recall a standard terminology. A discrete subgroup $\Gamma$ of $\mathrm{SL}(2,{\mathbb C})$ is called 
quasi-Fuchsian if the limit set of the action of $\Gamma$ on ${\mathbb C}{\mathbb P}^1$ is actually contained 
in a $\Gamma$--invariant Jordan curve.

\begin{lemma}\label{The1}
Let $\rho\,\colon\,\pi_1(X,\, x_0)\,\longrightarrow\, \mathrm{SL}(2,\C)$ be a discrete embedding. Then, the map
\[{\mathbb H}_\rho\,\colon\, \mathcal T(X)\,\longrightarrow\, {\mathcal H}_g\]
in \eqref{eP2} cannot be holomorphic, whenever $\rho$ happens to be quasi-Fuchsian.
\end{lemma}

\begin{proof}
Consider the smooth energy functional
\begin{equation}\label{ern}
\mathcal E_\rho\ \colon\ \mathcal T (X)\ \longrightarrow\ \R^{\geq0}
\end{equation}
which assigns to a given Riemann surface $\Sigma$ the (Dirichlet) energy (on the fundamental domain $ \Sigma$) of the
equivariant harmonic map with monodromy $\rho.$ Let $(E,\,\Phi)$ be the Higgs
bundle on $\Sigma$ corresponding to $\rho$; in other words, this map that sends
any $\Sigma$ to the above Higgs bundle $(E,\,\Phi)$ on $\Sigma$ gives ${\mathbb H}_\rho(\Sigma)$ (see
\eqref{eP2}). Then, the energy is given by
 \[\mathcal E_\rho(\Sigma)\,=\,\sqrt{-1}\,\int_\Sigma\text{tr}(\Phi\wedge\Phi^*).\]
 Moreover,
\[\det\Phi\,\in\, H^0(\Sigma,\,K_\Sigma^2)\,=\,T^*_\Sigma \mathcal T (X)\]
is a complex cotangent vector to ${\mathcal T}(X)$ at $\Sigma$,
and the derivative of the energy is given by
\begin{equation}\label{dE}
d_\Sigma\mathcal E_\rho\,\,=\,\,-\,4\cdot \Re\int_\Sigma(\text{tr}(\Phi^2),\, .)\,=\,
8\cdot \Re\int_\Sigma(\det(\Phi),\, .).
\end{equation}
In fact, the proof of \cite[Lemma 3.2]{Wolf} actually works verbatim in our situation. \\

It is shown in \cite[Section 8]{GW} that a discrete embedding $\rho\,\colon\,\pi_1(X,\,x_0)\,
\longrightarrow\, \mathrm{SL}(2,\mathbb C)$ is quasi-Fuchsian if and only if the function
$\mathcal E_\rho$ in \eqref{ern} is
proper. Assume that $\rho$ is quasi-Fuchsian. Since
$\mathcal E_\rho$ is proper, and it is bounded from below by $0,$ we conclude that
$\mathcal E_\rho$ is actually non-constant and attains its absolute minimum on $\mathcal T(X)$. \\

Assume that the map
$$\Sigma\,\longmapsto\,{\mathbb H}_\rho(\Sigma)\,=\,(E_\rho(\Sigma),\,\Phi_\rho(\Sigma))$$
in \eqref{eP2} is holomorphic. This implies that the map
\[\mathcal Q_\rho\,\colon \, \mathcal T (X)\,\longrightarrow\, \phi_* K^{\otimes 2}_\phi, \,\ \
\Sigma\,\longmapsto\, \det\Phi_\rho(\Sigma)\]
is holomorphic, where $K_\phi$ is the relative canonical line bundle for the projection $\phi$
in \eqref{e6}. Therefore,
\[\Sigma\,\,\longmapsto\, \,\int_\Sigma (\det\Phi_\rho(\Sigma)\wedge .)\]
is a holomorphic one-form on the Teichm\"uller space,
and consequently from \eqref{dE} we obtain that the energy functional $\mathcal E_\rho$ is given by the real
part of a holomorphic function (on the simply connected Teichm\"uller space). Thus, by the
maximum principle, $\mathcal E_\rho$ cannot attain its absolute minimum in $\mathcal T(X)$ contradicting that $\rho$ is quasi-Fuchsian. 
Therefore, ${\mathbb H}_\rho$ cannot be holomorphic.
\end{proof}

Next, we will extend Lemma \ref{The1} to the general case. For that purpose, we need a recent result by 
To\v{s}i\'{c} \cite{Tos} (see also Theorem \ref{thm:tos2} below) showing the energy $\mathcal E_\rho$ is 
strictly plurisubharmonic at points $\Sigma\,\in\,\mathcal T$ at which the spectral curve $S$ of the Higgs 
data ${\mathbb H}_\rho(\Sigma) \,=\, (\Sigma,\, (E,\,\Phi))$ is smooth (see \eqref{def:S} below for the 
definition of the spectral curve). Before proving Theorem \ref{The2} we first give a brief summary of the 
arguments and notions of \cite{Tos}.\\

\subsection{The kernel of the Levi form of the energy}\label{sec:tos}

The main input from \cite{Tos} is the computation of the kernel of the Levi form of the energy functional 
$${\mathcal E}\,=\,{\mathcal E}_\rho\,\colon\, {\mathcal T}\,\longrightarrow \,\R^{\geq 0}$$ associated to an 
irreducible representation $\rho\,\colon\, \pi_1(X,\, x_0) \, \longrightarrow\, \mathrm{SL}(2,\C).$ Our aim 
is to show that $\mathcal E_\rho$ cannot be pluriharmonic unless $\rho$ is unitary (in which case $\mathcal E_\rho$ is 
constant). Then, the Higgs data ${\mathbb H}_\rho(\Sigma)$ (see \eqref{eP2}) cannot depend holomorphically on 
$\Sigma\in\mathcal T(X).$\\
 
Consider the complex Hessian, i.e., the Levi form 
\[H(\mathcal E_\rho)\ :=\ d'd'' \mathcal E_\rho\ \in\ \Omega^{(1,1)}(\mathcal T(X))\]
of the energy functional with respect to the K\"ahler structure on the Teichm\"uller
space. Note that $H(\mathcal E_\rho)$ vanishes identically if and only if
$\mathcal E_\rho$ is pluriharmonic; see, for example, \cite{Dem}.

\begin{theorem}[{\cite[p.~7821, Theorem 1.6]{Tos}}]\label{cri1}
Let $\Sigma$ be a compact Riemann surface of genus $g\,\geq\, 2$ with underlying smooth surface $X$.
Let $\rho\,\colon\, \pi_1(X,\,x_0)\,\longrightarrow\,\mathrm{SL}(2,\mathbb C)$ be an irreducible representation
with corresponding Higgs pair $(\overline\partial,\,\Phi)$ at $\Sigma\,\in\,{\mathcal T}(X)$.
Then, the kernel $\mathcal K_\rho(\Sigma)$ of the Levi-form $H(\mathcal E_\rho)$ of $\mathcal E_\rho$ at $\Sigma$ is spanned by tangent vectors
$$\mu\,\in\, T_\Sigma{\mathcal T}(X) \,=\,H^1(\Sigma,\, T\Sigma)$$ such that there exists
$\xi\,\in\,\Gamma(\Sigma,\,\mathfrak{sl}(2,\mathbb C))$ with
\begin{equation}\label{kernelmuxi}
\overline{\partial} \xi\,=\,\mu\Phi\ \ \,\, \text{ and }\ \ \,\, [\xi,\,\Phi]\,=\,0.
\end{equation}
\end{theorem}

The proof of Theorem \ref{cri1} is a consequence of the Bochner technique combined with a
precise analysis of Toledo's proof \cite{Toledo} of the subharmonicity of $\mathcal E_\rho$.

 For an $\mathrm{SL}(2,\mathbb C)$ Higgs bundle $(\overline\partial,\, \Phi)$, the spectral curve is defined to be 
 \begin{equation}\label{def:S}S\,:=\,\{\omega\,\in\, K_\Sigma\,\big\vert\,\, -\omega^2\,=\,\det (\Phi)\}.\end{equation}
If $\det(\Phi)$ has only simple zeros, then $S$ is a smooth Riemann surface which
gives a double covering $\pi\,\colon\,S\, \longrightarrow\,\Sigma.$ Moreover,
$S$ naturally parametrizes the eigenvalues of $\Phi$ via $\omega$. To\v{s}i\'{c}
observed that the kernel $\mathcal K_\rho(\Sigma)$ of the Levi-form $H(\mathcal E_\rho)$ depends only on the
spectral data $(S,\,\omega),$ and not on the eigenlines of $\Phi$
(see \cite[p.~7821, Theorem 1.2]{Tos}). In the case at hand, this yields the following:

\begin{theorem}[{\cite{Tos}}]\label{thm:tos2}
Let $\Sigma$ be a compact Riemann surface of genus $g\,\geq\,2$ with underlying smooth surface $X$.
Let $\rho\,\colon \,\pi_1(X)\,\longrightarrow\,\mathrm{SL}(2,\mathbb C)$ be an irreducible representation
with corresponding Higgs pair $(\overline{\partial},\,\Phi)$ at $\Sigma\,\in\,\mathcal T(X)$. If the
spectral curve $S$ of $(\overline{\partial},\,\Phi)$ is smooth, then the energy is strictly pluri-sub-harmonic (in particular,
it is \textit{not} pluriharmonic).
\end{theorem}

\begin{proof}
We give To\v{s}i\'{c}'s proof for the $\mathrm{SL}(2,\C)$-case at hand: In the case of a Hitchin representation,
i.e., a Fuchsian representation, the corresponding statement about strictly pluri-sub-harmonicity was already
shown by Slegers in \cite[Proof of Theorem 1.1]{Sle}. 

Recall that for a given Riemann surface $\Sigma$, and any holomorphic quadratic differential $q\,\in\, 
H^0(\Sigma,\, K_\Sigma^2)$, there is a unique Hitchin representation on $\Sigma$ with Higgs field
$\Phi$, such that $\det (\Phi)\,=\, q$. In particular, this holds when $q$ has simple zeros and the
corresponding spectral curve $S\,=\,S(q)$ is smooth.

For a general 
representation, an element $\mu$ in the kernel of the Levi form gives rise to 
$\xi\,\in\,\Gamma(\Sigma,\,\mathfrak{sl}(2,\mathbb C))$ satisfying \eqref{kernelmuxi}. But $\xi$ determines (and is 
determined by) an odd function $x\, \colon\, \widehat \Sigma\,\longrightarrow\,\C$, given by the eigenvalues of
$\xi$ (the eigenlines 
of $\xi$ must be the same as the eigenlines of $\Phi$ by the second equation of \eqref{kernelmuxi}), and satisfy
\[(\pi^*\mu) \omega\,=\,\overline{\partial} x.\]
Therefore, the condition that $\mu$ is in the kernel of the Levi form depends only on the spectral data 
$(S,\,\omega)$ of $(\overline\partial,\, \Phi)$ (provided $S$ is smooth). Hence the result follows from the 
statement for the corresponding Hitchin representation with $q\,=\,\det\Phi$.
\end{proof}

\subsection{Elementary representations}\label{se-er}

We will also make use of the notion of {\em elementary representations}. In the case of $G\,=\,
\mathrm{SL}(2,\C)$, a representation is called elementary if and only if it is either unitary,
or reducible or it maps into a subgroup conjugated to the group generated by the embedding
\begin{equation}\label{eb}
\C^*\,\xhookrightarrow{}\, \mathrm{SL}(2,\C),\ \,\mu\,\longmapsto\,
\begin{pmatrix} \mu&0\\0&\mu^{-1}\end{pmatrix}
\end{equation}
and 
\[\begin{pmatrix} 0&-1\\1&0\end{pmatrix}\]
see \cite{WW} or \cite{GKM}. A representation which is not elementary is called non-elementary. 
A famous theorem of Gallo, Kapovich and Marden, \cite{GKM}, states that a surface group representation
into $\mathrm{SL}(2,\C)$ is non-elementary if and only if it is the (lifted) monodromy representation of
some projective structure on $X$. Additionally, it
 is proven by Wentworth and Wolf in \cite{WW} that an $\mathrm{SL}(2,\mathbb C)$-representation 
is non-elementary if and only if there is a Riemann surface $\Sigma\in\mathcal T(X)$ 
such that $\mathbb H_\rho(\Sigma)$ (see \eqref{eP2}) has a smooth spectral curve $S$.
It is exactly this characterization of non-elementary representations which will be used in the
proof of Theorem \ref{The2}.

\subsection{Non-existence for $G=\mathrm{SL}(2,\mathbb C)$}

Now we are in a position to prove the main theorem of this section.

\begin{theorem}\label{The2}
Let X be a surface of genus $g \,\geq\, 2$ and
$\rho\,\colon\,\pi_1(X,\, x_0)\,\longrightarrow\, \mathrm{SL}(2,\C)$ a
completely reducible representation. Then, the map
\[{\mathbb H}_\rho\,\colon\, \mathcal T (X)\,\longrightarrow\, {\mathcal H}_g\]
in \eqref{eP2} is holomorphic if and only if $\rho$ is unitary.
\end{theorem}

\begin{proof}
As in the proof of Lemma \ref{The1}, let $\rho \,\colon\, \pi_1(X,\, x_0)\,\longrightarrow\,
\mathrm{SL}(2,\C)$ be a completely reducible representation such that
\begin{itemize}
\item $\rho$ is \textit{not} unitary, and

\item the section ${\mathbb H}_\rho$ is holomorphic.
\end{itemize}

The second condition implies that $\mathcal 
E_\rho$ in \eqref{ern} must be pluriharmonic, as it is the real part of a holomorphic
function (compare with the proof of Lemma \ref{The1}).\\

By a result of Toledo \cite{Toledo}, the function $\mathcal E_\rho$ is always 
pluri-sub-harmonic. Moreover, by \cite[p.~7821, Theorem 2.1]{Tos}, see also Theorem \ref{thm:tos2} above, for
 every representation $\rho$ such that its Higgs data 
define a smooth spectral curve $S\,\longrightarrow\,\Sigma$ at some $\Sigma\,\in\,{\mathcal T}(X)$, the
function ${\mathcal E}_\rho$ is 
strictly pluri-sub-harmonic at $\Sigma$. Hence ${\mathcal E}_\rho$ is not pluriharmonic in that case. \\

Thus, the representations $\rho$ with holomorphic map ${\mathcal T}(X)\,\longrightarrow\, {\mathcal H}_G$ 
defined by $\Sigma\,\longmapsto\, \mathbb H_\rho(\Sigma)$ (see \eqref{eP2})
must give rise to a singular spectral curve at 
every $\Sigma\,\in\,{\mathcal T}(X)$. This property of a homomorphism --- that the corresponding spectral 
curve is singular for every $\Sigma\,\in\,{\mathcal T}(X)$ --- characterizes the elementary representations into 
$\mathrm{SL}(2,\mathbb C)$. In fact, this is the main result of \cite{WW}. Therefore, every 
$\mathrm{SL}(2,\mathbb C)$-representation $\rho$ with holomorphic map ${\mathbb H}_\rho$ must be an 
elementary representation.\\

Recall from Section \ref{se-er} the classification of elementary representations
into three types. To prove the theorem it suffices to exclude the following two cases:
\begin{enumerate}
\item $\rho$ is reducible.

\item The image of $\rho$ is contained in the subgroup generated by the diagonal matrices
and $\begin{pmatrix} 0&-1\\ 1&0\end{pmatrix}$.
\end{enumerate}

A reducible and totally reducible $\mathrm{SL}(2,\C)$--representation is conjugate to a diagonal
representation and it follows directly from \cite[p.~7821, Proposition 1.3]{Tos} that the energy
is not pluriharmonic (at every point $\Sigma \,\in\,\mathcal T(X)$ in Teichm\"uller space).\\

Therefore, it remains to show that the image of $\rho$ is not contained in the subgroup generated by the 
diagonal matrices and $\begin{pmatrix} 0&-1\\ 1&0\end{pmatrix}$.
Assume the converse, i.e., that the image of $\rho$ is contained in the subgroup generated by the 
diagonal matrices and $\begin{pmatrix} 0&-1\\ 1&0\end{pmatrix}$.
Observe that $\rho$ induces a ${\mathbb Z}/2{\mathbb Z}$-representation
\[\underline\rho\,\colon\, \pi_1(X,\,x_0)\,\longrightarrow\, {\mathbb Z}/2{\mathbb Z}\]
defined by
\[\underline\rho(\gamma)\,=\,\begin{cases} 0\, &\text{ if } \rho(\gamma) \text{ is diagonal }\\
1\, &\text{ if } \rho(\gamma) \text{ is off-diagonal.}\end{cases}\]
Let $K$ denote the kernel of $\underline \rho$, and let $\widetilde X$ be the universal covering of $X.$
Then $\widetilde X/K$ defines an unbranched 2-fold covering
\[\pi\,\colon\,\widehat{X}\,\longrightarrow \,X,\] if $\underline\rho$ is non-trivial. The compact
surface $\widehat X$ is connected and has genus $g_{\widehat X}\,=\,2\cdot \text{genus}(X)-1 \,=\, 2g-1$. By 
the construction of $\pi$, the induced representation $\pi^*\rho$ on $\widehat X$ is diagonal. The 
Teichm\"uller space of $\widehat X$ has complex dimension $6g-6$, while the Teichm\"uller space for $X$ is of 
dimension $3g-3.$ Note that any complex structure $J$ on $X$ gives rise to a unique complex
structure $\widehat J$ on $\widehat X$ such that the double covering $\pi\colon \widehat X\,
\longrightarrow\, X$ becomes holomorphic with respect to $\widehat J$ and $J$. Let $\widehat\Sigma
\,=\, (\widehat X, \widehat J)$ denote the Riemann surface corresponding to $\Sigma$ and let
\begin{equation}\label{ei}
\iota\,\,\colon\,\,{\mathcal T}(X)\,\longrightarrow\, {\mathcal T}({\widehat X}),
\ \ \, \Sigma\, \longmapsto\,\widehat\Sigma
\end{equation}
denote the corresponding map. Any element in the image of the map $\iota$ in
\eqref{ei} will be called a \textit{symmetric} Riemann surface; the reason for this terminology
is that we have a natural fixed-point-free involution $\sigma$ of $\widehat \Sigma$
satisfying the condition
\begin{equation}\label{eq:pisgn}
\pi\circ\sigma\ =\ \pi.
\end{equation}

The diagonal representation $\pi^*\rho$ gives rise to a diagonal Higgs pair on $\widehat \Sigma$
determined by a holomorphic 1-form 
\begin{equation}\label{eq:whomega}
\widehat\omega\,\in \,H^0(\widehat \Sigma,\, K_{\widehat\Sigma})
\end{equation}
and a flat unitary line bundle, unique up to sign. Note that
\[
0 \neq-\widehat\omega^2\,\,=\,\,\pi^*\det\Phi\,\,\in\,\, H^0(\widehat\Sigma,\,K_{\widehat\Sigma}^2).
\]
We consider the function $\mathcal E_{\pi^*\rho}\,\colon\,\mathcal T(\widehat X)\,\longrightarrow\, \R.$ This 
function is just twice the energy function of the corresponding $\mathbb C^*$-representation, and restricted 
to the space of symmetric Riemann surfaces, it is twice the energy $\mathcal E_{\rho}.$\\

The tangent space at $\widehat\Sigma$ of the space of symmetric Riemann surfaces is given by
\[\{\mu\,\in\, H^1(\widehat \Sigma,\, T\widehat\Sigma)
\,=\,T_{\widehat\Sigma}\mathcal T(\widehat X)\,\,\big\vert\,\, \sigma^*\mu\,=\,\mu\}.\]
If $\mathcal E_\rho$ would be pluriharmonic, every tangent vector to the subspace of $\sigma$-symmetric
Riemann surfaces in $\mathcal T(\widehat X)$
would lie in the kernel $\mathcal K(\widehat\Sigma)$ of the Levi form of $\mathcal E_{\pi^*\rho}$. From
the above discussion together with \cite[p.~7821, Proposition 1.3]{Tos}, that kernel is given by
\[\mathcal K_{\pi^*\rho}(\widehat\Sigma)\,=\,\{\mu\,\in\, H^1(\widehat \Sigma, \, T\widehat\Sigma)
\,=\,T_{\widehat{\Sigma}} \mathcal T(\widehat X)\,\big\vert\,\,\mu\,\widehat\omega \,=\,0\,
\in\, H^{0,1}(\widehat{\Sigma})\},\]
where $\widehat\omega$ is as in \eqref{eq:whomega}.\\

Consider $\underline\mu \,\in\, T_{\Sigma} {\mathcal T}(X)\,\cong\, H^0(\Sigma,\, K^2_\Sigma)^*$ such that
\[0\,\neq\,\int_\Sigma\underline\mu\det\Phi.\]
Then $\pi^*\underline\mu$ satisfies $$\sigma^*\pi^*\underline\mu\ =\ \pi^*\underline\mu$$ by
\eqref{eq:pisgn}, i.e., it represents a tangent vector to the space of symmetric Riemann
surfaces. On the other hand we obtain that
\[0\,\neq\,\int_\Sigma\underline\mu\det\Phi\,=\,-
\tfrac{1}{2}\int_{\widehat X}(\pi^*\underline\mu\,\widehat\omega)\widehat\omega.\]
Therefore, $\mathcal E_{\pi^*\rho}$ is not pluriharmonic along the complex submanifold of
symmetric Riemann surfaces. Since
\[\mathcal E_{\pi^*\rho}\circ \iota\,=\,2\mathcal E_{\rho},\]
this shows that $\mathcal E_{\rho}$ is also not pluriharmonic, and the proof is
complete.
\end{proof}

\section{Counter-examples in higher ranks}\label{sec:hr}

For the construction of counter-examples to Theorem \ref{The2} for higher rank flat connections, we
first recall a result by Koberda and Santharoubane \cite{KR}. They show that there are triples
\begin{equation}\label{et}
(X,\ N,\ \widehat{\rho})
\end{equation}
with the following properties:
\begin{itemize}
\item $X$ is a compact connected $C^\infty$ oriented surface of genus at least two,

\item $N$ is a (large) positive integer, and

\item $\widehat{\rho}$ is a homomorphism of groups
\begin{equation}\label{e3b}
\widehat{\rho}\, :\, {\rm MCG}^1(X) \, \longrightarrow\, \text{PGL}(N, {\mathbb C})
\end{equation}
\end{itemize}
(see \eqref{e7} for the definition of ${\rm MCG}^1(X)$) satisfying the following two conditions:
\begin{enumerate}
\item[(a)] the restriction of $\widehat{\rho}$ to the subgroup $\pi_1(X,\, x_0) \, \subset\, {\rm MCG}^1(X)$
(see \eqref{e8}) is irreducible, and

\item[(b)] there is an element $\gamma\, \in\, \pi_1(X,\, x_0)\, \subset\, {\rm MCG}^1(X)$ (see \eqref{e8}) such
that the adjoint action of $\widehat{\rho}(\gamma)\, \in\,
\text{PGL}(N, {\mathbb C})$ on the Lie algebra $\text{Lie}(\text{PGL}(N,
{\mathbb C}))$ has an eigenvalue $\lambda\, \in\, {\mathbb C}$ with
\begin{equation}\label{e4}
\big\vert\lambda\big\vert \ >\ 1.
\end{equation}
\end{enumerate}

\begin{remark}\label{rem1}
See \cite[p.~280, Theorem 3.1]{KR} and its proof for the above inequality
in \eqref{e4}. If the parameter $p$ in \cite[p.~280, Theorem 
3.1]{KR} is a prime number, then it was shown in \cite{KR} that $\widehat{\rho}$ is irreducible.
Irreducibility of $\widehat{\rho}$ was also proved in \cite[p.~26, Corollary A.4]{Go}.
\end{remark}

\begin{remark}
We note that there is no such $\widehat\rho$ for the case of $N\,=\, 2$; see \cite{BGMW}.
\end{remark}

Consider the restriction of $\widehat{\rho}$ to the subgroup $\pi_1(X,\, x_0)\, \subset\, {\rm MCG}^1(X)$
in \eqref{e8}; this restriction is denoted by
\begin{equation}\label{f1o}
\rho\ :\ \pi_1(X,\, x_0) \ \longrightarrow\ \text{PGL}(N, {\mathbb C}).
\end{equation}
Then \eqref{e4} immediately implies the following two:
\begin{enumerate}
\item The subgroup $\widehat{\rho}(\pi_1(X,\, x_0))\, \subset\, \text{PGL}(N, {\mathbb C})$ is infinite.

\item The homomorphism $\rho$ is \textit{not} unitary; in other words, for the natural inclusion map
$$\text{Hom}(\pi_1(X,\, x_0), \, \text{PU}(N))/\text{PU}(N)\,
\hookrightarrow\, {\mathcal R}_X(\text{PGL}(N,{\mathbb C}))$$
(see \eqref{e2}), we have
\begin{equation}\label{e5}
\rho\,\, \notin\,\, \text{Hom}(\pi_1(X,\, x_0), \, \text{SU}(N))/\text{SU}(N).
\end{equation}
\end{enumerate}

\begin{theorem}\label{mthm}
Let $X$ be a compact oriented surface of genus $g$, with $g\, \geq\, 4$.
There exists an integer $n$ (large), and an irreducible non-unitary representation $$\rho\ \colon
\ \pi_1(X,\,x_0)\ \longrightarrow\ \mathrm{SL}(n,\mathbb C)$$ such that the map
\[\mathbb H_\rho\,\,\colon\,\,\mathcal T(X)\,\,\longrightarrow\,\,\mathcal H_G\]
in \eqref{eP2} is holomorphic.
\end{theorem}

\begin{proof}
Take any triple as in \eqref{et} constructed in \cite{KR}
\begin{equation}\label{e7b}
(X,\ N,\ \widehat{\rho})
\end{equation}
such that $\text{genus}(X)\, \geq \, 4$. As in \eqref{f1o}
\begin{equation}\label{f1}
\rho\ :\ \pi_1(X,\, x_0) \ \longrightarrow\ \text{PGL}(N, {\mathbb C}).
\end{equation}
is the restriction of $\widehat{\rho}$ to the subgroup $\pi_1(X,\, x_0)\, \subset\, {\rm MCG}^1(X)$
in \eqref{e8}; so $\rho$ in \eqref{f1} satisfies the conditions (a) and (b) for the homomorphism
in \eqref{f1o}.

Let
\begin{equation}\label{e8b}
(E_\rho,\ \nabla^\rho)
\end{equation}
be the flat principal $\text{PGL}(N, {\mathbb C})$--bundle on $Z$
corresponding to $\rho$ in \eqref{f1} (see \eqref{e9}). Let 
\begin{equation}\label{st}
{\mathcal T}'(X)\,\, \subset\, \, {\mathcal T}(X)
\end{equation}
be the locus of Riemann surfaces that admit no nontrivial automorphism. The action of the
mapping class group ${\rm MCG}(X)$ on ${\mathcal T}(X)$ preserves ${\mathcal T}'(X)$, and,
moreover, the resulting action of ${\rm MCG}(X)$ on ${\mathcal T}'(X)$ is free. The quotient
\begin{equation}\label{q3}
{\mathcal M}_g\ :=\ {\mathcal T}'(X)/{\rm MCG}(X)
\end{equation}
is the moduli space of all Riemann surfaces of genus $g\,=\, \text{genus}(X)$ that do not admit any
nontrivial holomorphic automorphism. From the given condition that $g\, \geq\, 4$ it follows that
${\mathcal M}_g$ coincides with the smooth locus of the moduli space of smooth curves of
genus $g$. Let
\begin{equation}\label{ur}
\beta\ :\ {\mathcal M}^1_g\ \longrightarrow\ {\mathcal M}_g
\end{equation}
be the universal Riemann surface over ${\mathcal M}_g$; we note that there is a unique universal
Riemann surface over ${\mathcal M}_g$ because the Riemann surfaces in ${\mathcal M}_g$ do not admit
any nontrivial automorphism.

Although we will not need it here, it can be shown that the exact sequence in \eqref{e8}
is the long exact sequence of homotopy groups for the fibration in \eqref{ur}.

Define
$$
{\mathcal Z} \,\, :=\,\, \phi^{-1}({\mathcal T}'(X)) \,\, \subset\,\, Z
$$
(see \eqref{e6} and \eqref{st}); so ${\mathcal Z}$ is the universal automorphism-free Riemann surface over
${\mathcal T}'(X)$. Let
\begin{equation}\label{f2}
\widetilde{\phi}\, :\, {\mathcal Z}\, \longrightarrow\, {\mathcal T}'(X)
\end{equation}
be the restriction of $\phi$ to ${\mathcal Z}$. 

We note that the complement
\begin{equation}\label{co}
{\mathcal T}(X)\setminus {\mathcal T}'(X)\ \subset\ {\mathcal T}(X)
\end{equation}
is a complex analytic subspace of complex codimension at least two, because 
$g\, \geq\, 4$. Therefore, the inclusion map ${\mathcal T}'(X)\, \hookrightarrow\,
{\mathcal T}(X)$ induces an isomorphism of fundamental groups. Now, because ${\mathcal T}(X)$
is simply connected (recall that it is contractible), it follows that ${\mathcal T}'(X)$ is also
simply connected.

Consider the quotient map $\eta$ in \eqref{uc0}. Denote
\begin{equation}\label{uc}
\widehat{\mathcal Z} \ :=\ \eta^{-1}({\mathcal Z})\ \subset\ \eta^{-1}(Z)\ =\ \widehat{Z}.
\end{equation}
It was noted earlier that the map $\widehat\phi$ in \eqref{e6a} is equivariant (see \eqref{f3});
also, the action of ${\rm MCG}(X)$ on ${\mathcal T}(X)$ preserves ${\mathcal T}'(X)$.
Using these, it is deduced that the action of
${\rm MCG}^1(X)$ on $\widehat{Z}$ preserves $\widehat{\mathcal Z}$. Since $\widehat{\mathcal Z}/\pi_1(X,\, x_0)
\,=\, \mathcal Z$, from \eqref{e8} it follows that
\begin{equation}\label{q4}
\widehat{\mathcal Z}/{\rm MCG}^1(X)\ =\ {\mathcal Z}/{\rm MCG}(X)\ =\ {\mathcal M}^1_g
\end{equation}
(see \eqref{ur}).

It was noted earlier that the complement in \eqref{co} is a complex analytic subspace of complex codimension at least
two. From this it follows immediately that the complement $\widehat{Z}\setminus \widehat{\mathcal Z}\, \subset\,
\widehat Z$ is a complex analytic subspace of complex codimension at least two (this codimension clearly
coincides with the codimension in \eqref{co}). Recall that $\widehat{Z}$ is a contractible complex
manifold. As the complement $\widehat{Z}\setminus \widehat{\mathcal Z}\, \subset\,
\widehat Z$ is a complex analytic subspace of codimension at least two, it follows that 
$\widehat{\mathcal Z}$ is actually simply connected. Therefore, the quotient map
\begin{equation}\label{ep3}
\eta'\ :\ \widehat{\mathcal Z}\ \longrightarrow\, {\mathcal Z}\ =\ \widehat{\mathcal Z}/\pi_1(X,\, x_0)
\end{equation}
is a universal covering, where $\eta'$ is the restriction of $\eta$ in \eqref{uc0}, with Galois
group $\pi_1(X,\, x_0)$.

Consider the flat principal $\text{PGL}(N, {\mathbb C})$--bundle $(E_\rho,\, \nabla^\rho)$ on $Z$
in \eqref{e8b}. We will restrict it to
${\mathcal Z}\, \subset\, Z$. As its construction will be used, we recall it. Consider the trivial
principal $\text{PGL}(N, {\mathbb C})$--bundle
\begin{equation}\label{tb}
\widehat{\mathcal Z}\times\text{PGL}(N, {\mathbb C})
\ \longrightarrow \ \widehat{\mathcal Z}
\end{equation}
(recall the $\eta'$ in \eqref{ep3} is a universal covering map with Galois group $\pi_1(X,\, x_0)$).
The subgroup $\pi_1(X,\, x_0)$ in \eqref{e8} acts on $\widehat{\mathcal Z}\times\text{PGL}(N, {\mathbb C})$
as follows: The action of any $\gamma\, \in\, \pi_1(X,\, x_0)$ sends any $(z,\, A)\, \in\,
\widehat{\mathcal Z}\times\text{PGL}(N, {\mathbb C})$ to $(z\gamma,\, \rho(\gamma^{-1})A)$, where $\rho$ is
the homomorphism in \eqref{f1}; as noted above, the group ${\rm MCG}^1(X)$ acts on $\widehat{\mathcal Z}$ via
holomorphic automorphisms, and hence the subgroup $\pi_1(X,\, x_0)$ in \eqref{e8} acts on $\widehat{\mathcal Z}$.
The corresponding quotient
\begin{equation}\label{q1}
(\widehat{\mathcal Z}\times\text{PGL}(N, {\mathbb C}))/\pi_1(X,\, x_0)\ \longrightarrow\
\widehat{\mathcal Z}/\pi_1(X,\, x_0)\ =\ \mathcal Z
\end{equation}
is the restriction, to $\mathcal Z$, of the flat principal $\text{PGL}(N, {\mathbb C})$--bundle $(E_\rho,\,
\nabla^\rho)$ on $Z$ in \eqref{e8b}. The flat connection on it is given by the trivial connection
on the trivial bundle in \eqref{tb}. 
The flat principal $\text{PGL}(N, {\mathbb C})$--bundle in \eqref{q1} will be denoted by
\begin{equation}\label{fb1}
(E^1_\rho,\ \nabla^\rho_1),
\end{equation}
so $(E^1_\rho,\, \nabla^\rho_1)\ =\ (E_\rho,\ \nabla^\rho)\big\vert_{\mathcal Z}$.

The above action of the subgroup $\pi_1(X,\, x_0)$ on the trivial 
principal $\text{PGL}(N, {\mathbb C})$--bundle in \eqref{tb} extends to an action of the bigger group
${\rm MCG}^1(X)$ (see \eqref{e8}). The action of any $\gamma\, \in\, {\rm MCG}^1(X)$ sends any $(z,\, A)\,
\in\, \widehat{\mathcal Z}\times\text{PGL}(N, {\mathbb C})$ to $(z\gamma,\, \widehat{\rho}(\gamma^{-1})A)$, where
$\widehat\rho$ is the homomorphism in \eqref{e7b}; as noted before, the group ${\rm MCG}^1(X)$ acts on
$\widehat{\mathcal Z}$ via holomorphic automorphisms. The corresponding quotient
\begin{equation}\label{q2}
(\widehat{\mathcal Z}\times\text{PGL}(N, {\mathbb C}))/{\rm MCG}^1(X)\ \longrightarrow\
\widehat{\mathcal Z}/{\rm MCG}^1(X)\ =\ {\mathcal M}^1_g
\end{equation}
is a flat principal $\text{PGL}(N, {\mathbb C})$--bundle on ${\mathcal M}^1_g$ (see \eqref{q4}).
The flat principal $\text{PGL}(N, {\mathbb C})$--bundle in 
\eqref{q2} will be denoted by
\begin{equation}\label{fb2}
(E^2_\rho,\ \nabla^\rho_2).
\end{equation}

Let
\begin{equation}\label{ph}
\mathbf{q} \ :\ {\mathcal T}'(X) \ \longrightarrow\ {\mathcal M}_g
\ :=\ {\mathcal T}'(X)/{\rm MCG}(X)
\end{equation}
be the quotient map (see \eqref{q3}).

The following diagram of maps is Cartesian
\begin{equation}\label{cd}
\begin{matrix}
{\mathcal Z} & \xrightarrow{\,\,\,\Psi\,\,\,}& {\mathcal M}^1_g\\
\,\,\,\, \Big\downarrow \widetilde{\phi} && \,\,\, \Big\downarrow\beta\\
{\mathcal T}'(X) & \xrightarrow{\,\,\, {\mathbf q}\,\,\,} & {\mathcal M}_g
\end{matrix}
\end{equation}
where $\widetilde{\phi}$, $\mathbf{q}$ and $\beta$ are constructed in \eqref{f2},
\eqref{ph} and \eqref{ur} respectively. For any point $z\, \in\,{\mathcal T}'(X)$, the two Riemann surfaces
$\widetilde{\phi}^{-1}(z)$ and $\beta^{-1}({\mathbf q}(z))$ are isomorphic and they do not admit
any nontrivial holomorphic automorphism. So there is a unique holomorphic isomorphism between
$\widetilde{\phi}^{-1}(z)$ and $\beta^{-1}({\mathbf q}(z))$. These isomorphisms, as $z$ moves over
${\mathcal T}'(X)$, together produce the map $\Psi$ in \eqref{cd}.

We have
\begin{equation}\label{iso}
\Psi^* (E^2_\rho,\ \nabla^\rho_2)\, \ =\,\ (E^1_\rho,\ \nabla^\rho_1),
\end{equation}
(see \eqref{fb1}, \eqref{fb2}, \eqref{cd}). This follows from the fact that
$(E^2_\rho,\ \nabla^\rho_2)$ (respectively, $(E^1_\rho,\ \nabla^\rho_1)$) is
constructed using $\widehat{\rho}$ (respectively, $\rho$) and $\rho$ is a restriction of $\widehat{\rho}$.
More precisely, the identity map of $\widehat{\mathcal Z}\times\text{PGL}(N, {\mathbb C})$ descends to the
isomorphism in \eqref{iso} (see the constructions of $(E^2_\rho,\, \nabla^\rho_2)$ and
$(E^1_\rho,\, \nabla^\rho_1)$).

The adjoint action of $\text{PGL}(N, {\mathbb C})$ on its Lie algebra $\text{Lie}(\text{PGL}(N,
{\mathbb C}))$ embeds $\text{PGL}(N, {\mathbb C})$ into $\text{SL}(N^2-1, {\mathbb C})$. Let
\begin{equation}\label{f6}
({\mathcal E}_\rho,\ \widehat{\nabla}^\rho)
\end{equation}
be the flat $\text{SL}(N^2-1, {\mathbb C})$--bundle on ${\mathcal M}^1_g$ given by the flat 
$\text{PGL}(N, {\mathbb C})$--connection $(E^2_\rho,\, \nabla^\rho_2)$ in \eqref{fb2} using the
adjoint representation of
$\text{PGL}(N, {\mathbb C})$. We note that from the given condition that $\rho$ in 
\eqref{f1} is irreducible it follows that $({\mathcal E}_\rho,\, \widehat{\nabla}^\rho)$ is completely
reducible.\\

We use the adjoint representation to pass from $\text{PGL}(N, {\mathbb C})$--connections to
$\text{SL}(N, {\mathbb C})$--connections. If we wish to construct $\text{PGL}(N, {\mathbb C})$--connections
satisfying the condition in the theorem (instead of $\text{SL}(N, {\mathbb C})$--connections), then
there is no need for this use of adjoint representation.

Take a smooth compactification 
$\overline{\mathcal M}^1_g$ of ${\mathcal M}^1_g$. Thus, $\overline{\mathcal M}^1_g$ is an irreducible smooth
complex projective variety and ${\mathcal M}^1_g$ is a Zariski open subset of $\overline{\mathcal M}^1_g$
such that the complement
\begin{equation}\label{pd}
{\mathbb D}\ :=\ \overline{\mathcal M}^1_g\setminus {\mathcal M}^1_g\ \subset\ \overline{\mathcal M}^1_g
\end{equation}
is a simple normal crossing divisor of $\overline{\mathcal M}^1_g$.

Let $(\widetilde{\mathcal V},\, \widetilde{\Theta})$ be the parabolic
$\text{SL}(N^2-1, {\mathbb C})$--Higgs bundle on $\overline{\mathcal M}^1_g$ for the flat principal
$\text{SL}(N^2-1, {\mathbb C})$--bundle $({\mathcal E}_\rho,\ \widehat{\nabla}^\rho)$ on ${\mathcal M}^1_g$
in \eqref{f6}
(see \cite[Theorem 25.21]{Mo1}, \cite[p.~444, Theorem 5.17]{Mo2} and \cite[pp.~445--446, \S~5.3.3]{Mo2}). The
parabolic structure on the parabolic vector bundle $\widetilde{\mathcal V}$ is over the divisor
$\mathbb D$ (defined in \eqref{pd}). In particular, $\widetilde{\mathcal V}$ restricts to a usual algebraic
vector bundle ${\mathcal V}$ on ${\mathcal M}^1_g$, and the meromorphic Higgs field $\widetilde{\Theta}$ is
a regular Higgs field $\Theta$ on the vector bundle ${\mathcal V}$ over ${\mathcal M}^1_g$. In other words,
$({\mathcal V},\, \Theta)$ is a usual $\text{SL}(N^2-1, {\mathbb C})$--Higgs bundle over ${\mathcal M}^1_g$.

Consider the pulled back Higgs bundle 
$(\Psi^*{\mathcal V},\, \Phi^*\Theta)$ on $\mathcal Z$, where $\Psi$ is the map in \eqref{cd}. This 
gives a holomorphic family of Higgs bundles on the family of Riemann surfaces $\widetilde{\phi}\, :\, {\mathcal 
Z}\, \longrightarrow\, {\mathcal T}'(X)$ parametrized by ${\mathcal T}'(X)$ (see \eqref{f2}). From
\eqref{iso} it follows immediately that the corresponding flat $\text{SL}(N^2-1, {\mathbb C})$--bundle on
$\mathcal Z$ coincides with flat $\text{SL}(N^2-1, {\mathbb C})$--bundle given by
$(E^1_\rho,\ \nabla^\rho_1)$ in \eqref{fb1} (using the adjoint representation of $\text{PGL}(N,{\mathbb C})$).

Therefore, the flat $\text{SL}(N^2-1, {\mathbb C})$--bundle on $Z$ given by
$(E_\rho,\, \nabla^\rho)$ using the adjoint representation of $\text{PGL}(N, {\mathbb C})$
corresponds to a holomorphic family of $\text{SL}(N^2-1, {\mathbb C})$--Higgs bundles on
$Z$, because the subset $\mathcal Z\,\subset\, Z$ is open and dense in analytic topology. 

But by \eqref{e4} the so constructed family $(E_\rho,\, \nabla^\rho)$ of flat principal $\text{SL}(N^2-1,
{\mathbb C})$--bundles
on $Z$ using the adjoint representation of $\text{PGL}(N, {\mathbb C})$
is \textit{not} given by a family of flat principal $\text{SU}(N^2-1)$--bundles.
\end{proof}

\section*{Acknowledgements}

We thank the referee for helpful comments to improve the exposition.
We thank Takuro Mochizuki and Ramanujan Santharoubane for their help. The first-named
author is supported by a J. C. Bose Fellowship (JBR/2023/000003).


\end{document}